\title{Salajan's conjecture on discriminating terms in an exponential sequence}
\author{Pieter Moree and Ana Zumalac\'arregui}
\documentclass[12pt]{article}
\usepackage{amssymb, latexsym, amsfonts,amsmath}
\def\@ptsize{2}

\newtheorem{Thm}{Theorem}
\newtheorem{Con}{Conjecture}
\newtheorem{Prob}{Problem}
\newtheorem{Lem}{Lemma}

\newtheorem{Def}{Definition}
\newtheorem{cor}{Corollary}

\newtheorem{Prop}{Proposition}
\newcommand{\qed}{\hfill $\Box$}

\begin{document}
\date{}
\maketitle
{\def\thefootnote{}
}
\begin{abstract}
\noindent Given a sequence of distinct positive integers $v_1,v_2,\ldots$ and any positive integer $n$, the discriminator
$D_v(n)$ is defined as the smallest positive integer $m$ such $v_1,\ldots,v_n$ are pairwise
incongruent modulo $m$. We consider the discriminator for the
sequence $u_1,u_2,\ldots$, where $u_j$ equals the absolute value of 
$((-3)^j-5)/4$, that is $u_j=(3^j-5(-1)^j)/4$. We prove a 2012 conjecture of Sabin Salajan characterizing the discriminator
of the sequence $u_1,u_2,\ldots$. 
\end{abstract}
\section{Introduction}
Given a sequence of distinct positive integers $v_1,v_2,\ldots$ and any positive integer $n$, the discriminator
$D_v(n)$ is defined as the smallest positive integer $m$ such $v_1,\ldots,v_n$ are pairwise
incongruent modulo $m$. Browkin and Cao \cite{BC} relate it to cancellations algorithms similar to the
sieve of Eratosthenes.\\
\indent The main problem is to give an easy description or characterization of $D_v(n)$ (in many cases such a characterization does
not seem to exist).
Arnold, Benkoski and
McCabe \cite{ABM} might have
been the first to consider this type of problem (they introduced also the name). They considered the case
where $v_j=j^2$. Subsequently various authors, see e.g. \cite{BSW,PP,Sun,Zieve}, studied the discriminator for polynomial
sequences.\\
\indent It is a natural problem to study the discriminator for non-polynomial sequences. Very little work has been done in this direction. E.g. there are some conjectures due to Sun \cite{Sun}  in case 
$v_j=j!$, $v_j={\binom{2j}{j}}$ and $v_j=a^j$.\\ 
\indent In this paper we study the discriminator for a closely related sequence
$u_1,u_2,\ldots$ with $u_j=|(-3)^j-5|/4=(3^j-5(-1)^j)/4$.
This sequence satisfies the binary recurrence
$u_n=2 u_{n-1} + 3 u_{n-2}$, for every $n \geq 3$. with
starting values $u_1=2$ and $u_2=1$.
The first few terms are $$2, 1, 8, 19, 62, 181, 548, 1639, 4922,\ldots$$
Note that for $j\ge 2$ we have $u_{j+1}>u_j$ and that all $u_j$ are distinct.
It is almost immediate that the terms are of alternating parity.
Since all $u_j$ are distinct the number
$$D_S(n)=\min\{m\ge 1:u_1,\ldots,u_n{\rm~are~pairwise~distinct~modulo~}m\}$$
is well-defined. 
Note that $D_S(n)\ge n$. In Table 1 we give the values of $D_S(n)$ for $1\le n\le 32768$ (with the powers
of $5$ underlined).\\
\centerline{\bf TABLE 1}
\begin{center}
\begin{tabular}{|c|c|c|c|}
\hline
range & value & range & value \\
\hline
$1$ & $1$ & $129-256$ & $256$\\
\hline
$2$ & $2$ & $257-512$ & $512$\\
\hline
$3-4$ & $4$ & $513-1024$ & $1024$\\
\hline
$5-8$ & $8$ & $1025-2048$ & $2048$\\
\hline
$9-16$ & $16$ & $2049-2500$ & $\underline{3125}$\\
\hline
$17-20$ & $\underline{25}$ & $2501-4096$ & $4096$\\
\hline
$21-32$ & $32$ & $4097-8192$ & $8192$\\
\hline
$33-64$ & $64$ & $8193-12500$ & $\underline{15625}$\\
\hline
$65-100$ & $\underline{125}$ & $12501-16384$ & $16384$\\
\hline
$101-128$ & $128$ & $16385-32768$ & $32768$\\
\hline
\end{tabular}
\end{center}
Based on this table Sabin Salajan, 
 who at the time was
an intern with the first author, proposed a conjecture that we will
prove in this paper to be true. The first author had asked Sabin to
find second order linear recurrences for which the discriminator values have 
a nice structure. After an extensive search Sabin came up with the sequence
$u_1,u_2,\ldots$. For convenience we call this sequence the
{\it Salajan sequence} $S$ and its associated discriminator $D_S$ the {\it Salajan discriminator}.
If $m=D_S(n)$ for some $n\ge 1$, then we say that $m$ is a {\it Salajan value}, otherwise it
is a {\it Salajan non-value}.
\begin{Thm}
\label{main}
Let $n\ge 1$.
Put $e=\lceil{\log_2(n)}\rceil$ and $f=\lceil{\log_{5}({5n/4)}}\rceil$.
Then $$D_S(n)=\min\{2^e,5^f\}.$$
\end{Thm}
\begin{cor}
If the interval $[n,5n/4)$ contains a power of $2$, say $2^a$, then
we have $D_S(n)=2^a$.
\end{cor}
Note that $2^e$ is the smallest power of $2$ which is $\ge n$ and
that $5^f$ is the smallest power of $5$ which is $\ge 5n/4$.\\
\indent {}From Table 1 one sees that not all powers of $5$ are Salajan values.
Let ${\cal F}$ be the set of integers $b\ge 1$ such that the interval $[4\cdot5^{b-1},5^b]$ does not
contain a power of $2$. Then it is not difficult to show that the image of
$D_S$ is given by $\{2^a:a\ge 0\}\cup \{5^b:b\in {\cal F}\}$. Using Weyl's criterion one
can easily establish (see Section \ref{Izabela}) the following
proposition.
\begin {Prop}
\label{iza}
As $x$ tends to infinity we have
 $\#\{b\in {\cal F}:b\le x\}\sim \beta x$, with $\beta = 3-\log 5/\log 2 = 0.678\ldots$. 
\end{Prop}

\section{Strategy of the proof of Theorem \ref{main}}
For the benefit of the reader we describe the strategy of the (somewhat lengthy) proof
of Theorem \ref{main}.\\
\indent We first show that if $2^e\ge n$ and $5^f\ge 5n/4$, then
$D_S(n)\le \min\{2^e,5^f\}$. This gives us the absolutely crucial upper bound
$D_S(n)<2n$.\\
\indent Next we study the periodicity of the sequence modulo $d$ and determine
its period $\rho(d)$. The idea is to use the information so obtained to show that
many $d$ are Salajan non-values.
In case $3\nmid d$ the sequence turns out to be purely periodic with even
period that can be given precisely. This is enough for our purposes as
we can show that $3|D_S(n)$ does not occur. 

Now we restrict to the $d$ with $3\nmid d$.
Using that $D_S(n)<2n$ one easily sees that if $\rho(d)\le d/2$, then $d$ is a  
Salajan non-value. The basic property (\ref{per2}) of the period together with the evenness of the period now
excludes composite values of $d$. Thus we have $d=p^m$, with $p$ a prime.

In order for $\rho(p^m)>p^m/2$ to hold we find that we must have ord$_9(p)=(p-1)/2$, that
is $9$ must have maximal possible order modulo $p$. 
Moreover, $9$ must have maximal possible order modulo $p^m$, that is ord$_9(p)=\varphi(p^m)/2$.
(A square cannot have a multiplicative order
larger than $\varphi(p^m)/2$ modulo $p^m$.) 
This is about as far as the study of the periodicity will
get us. To get further we will use a more refined tool, the {\it incongruence index}. Given an integer $m$,
this is the maximum $k$ such that $u_1,\ldots,u_k$ are pairwise distinct modulo $m$. We
write $\iota(m)=k$. For $3\nmid m$, $\iota(m)\le \rho(m)$. Using that $D_S(n)<2n$ one 
notes that if $\iota(d)\le d/2$, then $d$ is a Salajan non-value. 

For the primes $p>3$ 
we show by a lifting argument that if $\iota(p)<\rho(p)$, 
then $p^2,p^3,\ldots$ are Salajan non-values.
Likewise, we prove that if $\iota(p)\le p/2$, then $p,p^2,p^3,\ldots$ are Salajan non-values.
We then show that except for $p=5$, all primes with ord$_9(p)=(p-1)/2$ satisfy $\iota(p)<\rho(p)$.
At this point we are left with the primes $p>5$ satisfying ord$_9(p)=(p-1)/2$ as only possible Salajan values.
Then using classical exponential sums techniques, and some combinatorial arguments, we infer that 
$\iota(p)<4p^{3/4}$. Using this bound, after some computational work, we then conclude that
$\iota(p)\le p/2$ for every $p>5$.

Thus we are left with $D_S(n)=2^a$ for some $a$ or
$D_S(n)=5^b$ for some $b$. By Lemma \ref{prop2n} and 
Lemma \ref{prop5n} it now follows that $2^a\ge n$ and $5^b\ge 5n/4$. This
then completes the proof.


\section{Preparations for the proof}
\label{preparations}

We will show that $2^e$ with $2^e\ge n$ and $5^f$ with $5^f\ge 5n/4$ are
admissible discriminators. That is, we will show that the sequence
$u_1,\dots, u_n$ lie in distinct residue classes modulo
$2^e$ and in distinct residue classes modulo $5^f$.

Let $p$ be a prime. If $p^a|n$ and $p^{a+1}\nmid n$, then we
put $\nu_p(n)=a$. The following result is well-known, for a proof see,
e.g., Beyl \cite{Beyl}.
\begin{Lem}
\label{Beyl}
Let $p$ be a prime, $r\ne -1$ an integer satisfying $r\equiv 1({\rm mod~}p)$ and $n$ a natural number.
Then
$$
\nu_p(r^n-1) =
\begin{cases}
\nu_2(n)+\nu_2(r^2-1)-1 & \text{if $p=2$ and $n$ is even};\\
\nu_p(n)+\nu_p(r-1) & \text{otherwise}.
\end{cases}
$$
\end{Lem}
\begin{cor}
\label{hakbijl}
Let $f\ge 2$ and $p$ be an odd prime.
If $g$ is a primitive root modulo $p$, then $g$ is 
a primitive root modulo $p^f$ if and only if $g^{p-1}\not\equiv 1({\rm mod~}p^2)$. 
\end{cor}
Corollary \ref{hakbijl} is a classical result from elementary number
theory. For an alternative proof see, e.g., Apostol \cite[Theorem 10.6]{Apostel}.
\begin{Prop}
\label{alpha}
Let $n\ge 1$ be an integer, $p$ a prime and put $e_p=\lfloor{\log_p(n-1)}\rfloor$.
Let $r\equiv 1({\rm mod~}p)$ be
an integer $\ne -1$. Put $r_p=\nu_p(r-1)$. If $p=2$, we assume in addition that $r$ is a square.
The integers $r,\ldots,r^{n}$ are pairwise distinct modulo $p^{e_p+r_p+1}$.
\end{Prop}
{\it Proof}. Write $m=p^{e_p+r_p+1}$. Let $1\le i<j\le n$ and suppose that $r^{i}\equiv r^{j}({\rm mod~}m)$, thus
$r^{j-i}\equiv 1({\rm mod~}m)$ and hence $\nu_p(r^{j-i}-1)\ge e_p+r_p+1$. Note that $\nu_p(k)\le e_p$ for $1\le k\le n-1$. Thus
$\nu_p(j-i)\le e_p$ and, by Lemma \ref{Beyl}, we deduce that $\nu_p(r^{j-i}-1)\le e_p+r_p$.
Contradiction.\qed
\begin{cor}\label{9and81}$~$\\
The integers $9,\ldots,9^{n}$ are pairwise distinct modulo $2^{e_2+4}$.\\
The integers $81,\ldots,81^n$ are pairwise distinct modulo $5^{e_5+2}$.
\end{cor}

\begin{Lem}
\label{prop2n}
Let $n \geq 2$ be an integer with $n\le 2^m$.  Then, we have that $u_1,\dots,u_n$ are pairwise distinct modulo $2^m$.
\end{Lem}
{\it Proof}. For $n=2$ the result is obvious. So assume that $n\ge 3$.
Since the terms of the sequence alternate between even and odd, it suffices to compare the remainders $({\rm mod~}2^m)$ of
the terms having an index with the same parity. Thus assume that we have
$$u_{2j+\alpha} \equiv u_{2k+\alpha}({\rm mod~}2^m){\rm ~with~}1\le 2j+\alpha<2k+\alpha \leq n,~\alpha\in \{1,2\}.$$ It follows from this that $9^{k-j}\equiv 1 ({\rm mod~}2^{m+2})$.
We have $\nu_2(9^{k-j}-1)=\nu_2(k-j)+3$ by Lemma \ref{Beyl}. Further, 
$2k-2j \leq n-1<2^m$, so $\nu_2(k-j) \leq m-2$ (here we used that $n\ge 3$).
Therefore  $\nu_2(9^{k-j}-1)=\nu_2(k-j)+3 \leq (m-2) +3 = m+1$, which implies that $9^{k-j}-1$ cannot be
divisible by $2^{m+2}$. Contradiction. \qed\\

\noindent {\tt Remark}. The incongruence of $u_i$ and $u_j$ (mod $2^m$) with $i$ and $j$ of the
same parity and $1\le i<j\le n$ is equivalent with $9,9^2,\ldots,9^{\lfloor (n-1)/2\rfloor}$ 
being pairwise incongruent mod $2^m$. Using this observation and Corollary \ref{9and81} we obtain 
an alternative proof of Lemma \ref{prop2n}.\\

On noting that trivially $D_S(n)\ge n$ and that  
for $n\ge 2$ the interval $[n,2n-1]$ always contains some power of $2$, we
obtain the following corollary to Lemma \ref{prop2n}.
\begin{cor}
\label{bertie}
We have $n\leq D_S(n)\le 2n-1$.
\end{cor}

\begin{Lem}
\label{prop5n}
The integers $u_1,\dots,u_n$ are pairwise distinct modulo $5^{m}$ iff $$5^m\ge 5n/4.$$
\end{Lem}
{\it Proof}. If $5^m<5n/4$, then $1+4\cdot 5^{m-1}\le n$. By Lemma \ref{Beyl} we have
$$81^{5^{m-1}}\equiv 1({\rm mod~}5^m),$$ which ensures that 
$u_1\equiv u_{1+4\cdot 5^{m-1}}({\rm mod~}5^m)$. Next let us assume that $5^m\ge 5n/4$. This
ensures that $m\ge 1$.
The remainders of the sequence modulo $5$ are $2,1,4,3,2,1,\dots$ and so the sequence has
period $4$ modulo $5$. Thus we may assume that $m\ge 2$.
It suffices to show that $u_{j_1}\not\equiv u_{k_1}({\rm mod~}5^{m})$ with $1\le j_1< k_1\le n$ in the same congruence class modulo $4$. We will argue by contradiction.
Thus we assume that
$$u_{4j+\alpha} \equiv u_{4k+\alpha}({\rm mod~}5^{m}){\rm ~with~}1\le 4j+\alpha<4k+\alpha \leq n,~\alpha\in \{1,2,3,4\}.$$ {}From this it follows that $81^{k-j}\equiv 1 ({\rm mod~}5^{m})$, where
$k-j\le (n-\alpha)/4<n/4\le 5^{m-1}$ by hypothesis and hence $\nu_5(k-j)\le m-2$.
On invoking Lemma \ref{Beyl} we now infer that $\nu_5(81^{k-j}-1)=\nu_5(k-j) + 1 \leq m-2+1=m-1$. Contradiction. \qed\\

\noindent {\tt Remark}. The incongruence of $u_i$ and $u_j$ (mod $5^m$) with $i$ and $j$ 
in the same residue class modulo $4$ and $1\le i<j\le n$ is equivalent with $81,81^2,\ldots,81^{\lfloor (n-1)/4\rfloor}$ 
being pairwise incongruent mod $5^m$. Using this observation and Corollary \ref{9and81} we obtain 
an alternative proof of Lemma \ref{prop5n}.\\

\noindent In order to determine whether a given $m$ discriminates $u_1,\ldots,u_n$ modulo $m$, we can separately
consider whether $u_i\not\equiv u_j({\rm mod~}m)$ with $1\le i<j\le n$ of the same parity (case 1)
and with distinct parity (case 2). The first case is easy and covered by Lemma \ref{oneventje}, the second case
is trivial in case $m$ is a power of $2$ or $5$, but in general much harder than the first case.
\begin{Lem}
\label{oneventje}
Suppose that $3\nmid m$ and $1\le \alpha\le n$. We have $u_i\not\equiv u_j({\rm mod~}m)$ for
every pair $(i,j)$ satisfying $\alpha\le i< j\le n$ with $i\equiv j({\rm mod~}2)$
iff ord$_9(4m)>(n-\alpha)/2$.
\end{Lem} 
{\it Proof}. We have $u_i\not\equiv u_{i+2k}({\rm mod~}m)$ iff $9^k\not\equiv 1({\rm mod~}4m)$. Thus
$u_i\not\equiv u_j({\rm mod~}m)$ for every pair $(i,j)$ with $\alpha\le i<j\le n$ and
$i\equiv j({\rm mod~}2)$ iff $9^k\not\equiv 1({\rm mod~}4m)$ for $1\le k\le (n-\alpha)/2$. \qed\\

\noindent {\it Alternative proof of Lemma} \ref{prop2n}. If $i$ and $j$ are of different parity, 
then $u_i\not\equiv u_j({\rm mod~}2)$. Hence we may assume that $i$ and $j$ are of
the same parity. On invoking Lemma \ref{oneventje} we then obtain that $u_1,\ldots,u_{n}$ are
distinct modulo $2^m$ iff ord$_9(2^{m+2})>(n-1)/2$. By Lemma \ref{Beyl} we have
ord$_9(2^{m+2})=2^{m-1}$, concluding the proof. \qed\\

\noindent {\it Alternative proof of Lemma} \ref{prop5n}. The remainders of the sequence 
modulo $5$ are $2,1,4,3,2,1,\dots$ and so terms $u_i$ and $u_j$ with 
$i$ and $j$ of different parity are incongruent. Now by Lemma \ref{oneventje}  the 
 integers $u_1,\dots,u_n$ are pairwise distinct modulo $5^{f}$ iff 
ord$_9(4\cdot 5^f)>(n-1)/2$. Since $3$ is a primitive root modulo $5$ 
and $3^4\not\equiv 1({\rm mod~}5^2)$, we have by Corollary \ref{hakbijl} that
$3$ is a primitive root modulo $5^f$ and hence ord$_3(5^f)=4\cdot 5^{f-1}=\varphi(5^f)$, with 
$\varphi$ Euler's totient function.
On making use of the trivial observation that, for integers $m$ coprime to $3$,
\begin{equation}
\label{lcmorder}
2{\rm ord}_9(4m)={\rm lcm}(2,{\rm ord}_3(4m)),
\end{equation}
we infer that ord$_9(4\cdot 5^f)={\rm ord}_9(5^f)={\rm ord}_3(5^f)/2=2\cdot 5^{f-1}$.
The proof is now finished by noting that the condition ord$_9(4\cdot 5^f)>(n-1)/2$ is
equivalent to $5^f\ge 5n/4$. \qed

\section{Periodicity and discriminators}
\subsection{Generalities}
We say that a sequence of integers $\{v_j\}_{j=1}^{\infty}$ is {\it (eventually) periodic} 
modulo $d$ if there exist integers $n_0\ge 1$ and $k\ge 1$ such that
\begin{equation}
\label{per1}
v_n\equiv v_{n+k}({\rm mod~}d)
\end{equation}
for every $n\ge n_0$. 
The minimal choice for $n_0$ is called the {\it pre-period}.
The smallest $k\ge 1$ for which (\ref{per1}) holds
for every $n\ge n_0$ is said to be the {\it period} and denoted by $\rho_v(d)$.
In case we can take $n_0=1$ we
say that the sequence is {\it purely periodic} modulo $d$.

\smallskip

 Let $\{v_j\}_{j=1}^{\infty}$ be a second order linear recurrence with the 
two starting values and the coefficients of the defining equation being integers.
Note that, for a given $d$, there must be a pair $(a,b)$ such hat $a\equiv v_n$ and $b\equiv v_{n+1}$ modulo $d$ 
for infinitely many $n$. Since a pair of consecutive terms determines uniquely all subsequent ones,
it follows that the sequence is periodic modulo $d$. If we consider $n$-tuples instead of pairs modulo $d$, we
see that an $n$th order linear recurrence with the 
$n$ starting values and the coefficients of the defining equation  being integers, is always periodic
modulo $d$.\\
\indent If a sequence $v$ is periodic modulo $d_1$ and modulo $d_2$ and
$(d_1,d_2)=1$, then we obviously have
\begin{equation}
\label{per2}
\rho_v(d_1d_2)={\rm lcm}(\rho_v(d_1),\rho_v(d_2)).
\end{equation}
If the sequence is purely periodic modulo $d_1$ and modulo $d_2$ and
$(d_1,d_2)=1$, then it is also purely periodic modulo $d_1d_2$.
Another trivial property of $\rho_v$ is that if 
the sequence $v$ is periodic modulo $d_2$, then for
every divisor $d_1$ of $d_2$ we have
\begin{equation}
\label{per3}
\rho_v(d_1)|\rho_v(d_2).
\end{equation}
\indent The following result links the period with the discriminator. Its moral is that
if $\rho_v(d)$ is small enough, we cannot expect $d$ to occur as $D_v$-value, i.e. $d$ does not 
belong to the image of $D_v$.

\begin{Lem}
\label{gee}
Assume that $D_v(n)\le g(n)$ for every $n\ge 1$ with $g$ non-decreasing.
Assume that the sequence $v$ is purely periodic modulo $d$ with period $\rho_v(d)$. If $g(\rho_v(d))<d$, then $d$ is a $D_v$-non-value.
\end{Lem}
{\it Proof}. Since $v_1\equiv v_{1+\rho_v(d)}({\rm mod~}d)$ we must have $\rho_v(d)\ge n$.
Suppose that $d$ is a $D_v$-value, that is for some $n$ we have $D_v(n)=d$. Then
$d=D_v(n)\le g(n)\le g(\rho_v(d))$. Contradiction. \qed

\subsection{Periodicity of the Salajan sequence}
\label{peri}
The purpose of this section is to establish Theorem~\ref{generalperiod}, which gives an explicit formula for the 
period $\rho(d)$ and the pre-period for the Salajan sequence.
Since it
is easy to show that $3\nmid D_S(n)$, it would be actually enough to study 
those integers $d$ with $3\nmid d$ (in which case the Salajan sequence is purely 
periodic modulo $d$). However, for completeness we discuss the periodicity of
the Salajan sequence for {\it every} d.
\begin{Thm}
\label{generalperiod}
Suppose that $d>1$. Write $d=3^{\alpha}\cdot \delta$ with $(\delta,3)=1$. 
The period of the Salajan sequence modulo $d$, $\rho(d)$, exists
and satisfies $\rho(d)=2{\rm ord}_9(4\delta)$. The pre-period equals
$\max(1,\alpha)$.
\end{Thm}
\begin{cor}
\label{9}
The Salajan sequence is purely periodic
iff $9\nmid d$. 
\end{cor}
\begin{Lem}
Write $d=3^{\alpha}\cdot \delta$ with $(\delta,3)=1$. The Salajan sequence is purely periodic
iff $9\nmid d$. Furthermore, if $9\nmid d$, then $\rho(d)\,|\,2{\rm ord}_9(\delta)$.
\end{Lem}
{\it Proof}. Since $u=2,{\overline{1,8}}({\rm mod~}9)$ the condition $9\nmid d$ is 
necessary for the Salajan sequence to be purely periodic modulo $d$. 

We will now show that it is also sufficient.
Let us first consider the case where $\alpha=0$. We note
that $u_n\equiv u_{n+2k} ({\rm mod~}d)$ iff $3^n\equiv 3^{n+2k}({\rm mod~}4\delta)$. It
follows that $\rho(d)\,|\,2{\rm ord}_9(4\delta)\,|\,2k$.
If $\alpha=1$, then we use (\ref{per2}) and the observation that 
$2=\rho(3)\,|\,2{\rm ord}_9(4\delta)$. \qed\\

\noindent {\tt Remark}. The above proof shows that if $\rho(d)$ is even, then $\rho(d)=2{\rm ord}_9(4\delta)$.

\begin{Lem}
\label{even}
Assume that $9\nmid d$ and $d>1$. The Salajan sequence is purely periodic with
period $\rho(d)=2{\rm ord}_9(4\delta)$, where $d=3^{\alpha}\cdot \delta$ with $(\delta,3)=1$.
\end{Lem} 
{\it Proof}. By the previous remark it suffices to show that $\rho(d)$ is even. If $\alpha=1$, then $2=\rho(3)\,|\,\rho(d)$
(here we use (\ref{per3})) and we are done, so we may assume that $\alpha=0$. If $5\,|\,d$, then
$4=\rho(5)\,|\,\rho(d)$ and so we may assume that $(5,d)=1$. Suppose that $\rho(d)$ is odd. Then
\begin{equation}
\label{per4}
u_{n}\equiv u_{n+\rho(d)}({\rm mod~}d)
\end{equation}
iff
$3^n-5(-1)^n\equiv 3^{n+\rho(d)}+5(-1)^n({\rm mod~}4d)$ iff $5^*(1-3^{\rho(d)})/2\equiv (-3)^{-n}({\rm mod~}2d)$,
where $5^*$ is the inverse of 5 modulo $2d$. Now if (\ref{per4}) is to hold for every $n\ge 1$, then
$(-3)^n$ assumes only one value as $n$ ranges 
over the positive integers. Since $(-3)^{\phi(2d)}\equiv 1({\rm mod~}2d)$ we must have
$(-3)^n\equiv 1({\rm mod~}2d)$ for every $n\ge 1$. 
This implies that $d=2$ or $d=1$.
Since $5^*(1-3^2)/2\not\equiv 1({\rm mod~}4)$ it follows that $d=1$. Contradiction. \qed\\

\noindent {\it Proof of Theorem} \ref{generalperiod}. It is an easy observation that
modulo $3^{\alpha}$ the Salajan sequence has pre-period $\max(\alpha,1)$ and period two.
This in combination with Lemma \ref{even} and (\ref{per2}) then completes the proof. \qed

\subsection{Comparison of $\rho(d)$ with $d$}
\begin{Lem}
\label{liftje}
Let $p>3$. We have $\rho(p^m)\,|\,\rho(p)p^{m-1}$.
\end{Lem} 
{\it Proof}. Since $3^{\rho(p)}\equiv 1({\rm mod~}p)$ we
have $3^{\rho(p)p^{m-1}}\equiv 1({\rm mod~}p^m)$ and, provided that $\rho(p)$ is even,
this implies that $u_k\equiv u_{k+\rho(p)p^{m-1}}({\rm mod~}p^m)$ for every $k\ge 1$. \qed

\begin{cor}
\label{liftje2}
Either $\rho(p^2)=\rho(p)$ or $\rho(p^2)=p\rho(p)$.
\end{cor}
{\it Proof}. We have $\rho(p)|\rho(p^2)|p\rho(p)$.

\begin{Lem}
\label{previous}
We have $\rho(2^e)=2^e$ and $\rho(3^e)=2$. If $p$ is odd, then $\rho(p^e)|\varphi(p^e)$. 
\end{Lem}
{\it Proof}. {}From Lemma \ref{Beyl} and Lemma \ref{even} we infer
that ord$_9(2^{e+2})=2^{e-1}$ and hence $\rho(2^e)=2^e$. 
For $n$ large enough modulo $3^e$ the sequence alternates between $-5/4$ and $5/4$ modulo $3^e$.
Since these are different residue classes, we have $\rho(3^e)=2$.\\
\indent It remains to prove the final claim. 
If $p=3$ it is clearly true and thus we may assume that $p>3$. 
Note that $\rho(p^e)=2{\rm ord}_9(4p^e)=2{\rm ord}_9(p^e)$ and that
$2{\rm ord}_9(p^e)\,|\,2(\varphi(p^e)/2)=\varphi(p^e)$. \qed

\begin{cor}
\label{previous1}
We have $\rho(d)\le d$.
\end{cor}

\begin{Lem}
\label{uppierho}
Suppose that $d_1,d_2>1$ and $(d_1,d_2)=1$. Then
$$\rho(d_1d_2)\le \rho(d_1)\rho(d_2)/2\le d_1d_2/2.$$
\end{Lem}
{\it Proof}. We have $\rho(d_1d_2)={\rm lcm}(\rho(d_1),\rho(d_2))$. 
By Lemma \ref{even} both $\rho(d_1)$ and $\rho(d_2)$ are even. It thus follows that
$\rho(d_1d_2)\le \rho(d_1)\rho(d_2)/2$. The final estimate follows 
by Corollary \ref{previous1}. \qed

\section{Non-values of $D_S(n)$}
\label{non-values}
Recall that if $m=D_S(n)$ for some $n\ge 1$ we call $m$ a Salajan value and
otherwise a Salajan non-value.\\
\indent Most of the following proofs rely on the simple fact that for certain
sets of integers we have that if $u_1,\dots,u_n$ are in 
$n$ distinct residue classes modulo $m$, then $m\ge 2n$ 
contradicting Corollary \ref{bertie}.

\subsection{$D_S(n)$ is not a multiple of $3$}
\begin{Lem}
\label{notdrie}
We have $3\nmid D_S(n)$.
\end{Lem}
{\it Proof}. We argue by contradiction and so assume that 
$D_S(n)=3^{\alpha}m$ with $(m,3)=1$ and $\alpha \ge 1$.
Since by definition $u_{\alpha}\not\equiv u_{\alpha+2t}({\rm mod~}3^{\alpha}m)$
for $t=1,\ldots,\lfloor{(n-\alpha)/2}\rfloor$ and  $u_{\alpha}\equiv u_{\alpha+2t}({\rm mod~}3^{\alpha})$ 
for every $t\ge 1$, it follows that $u_{i}\not\equiv u_{j}({\rm mod~}m)$ 
with $\alpha \le i<j\le n$ and $i$ and $j$ of the same parity. By Lemma \ref{oneventje}
it then follows that ord$_9(4m)>(n-\alpha)/2$. By Lemma \ref{even}, Corollary \ref{previous1}
and Corollary \ref{bertie} we then find that $n-\alpha+1\le 2{\rm ord}_9(4m)=\rho(m)\le m\le 2n/3^{\alpha}$. 
This implies
that $n\le 3^{\alpha}(\alpha-1)/(3^{\alpha}-1)$. On the other hand, by
 Corollary \ref{bertie} we have
$3^{\alpha}m\le 2n$ and hence $n\ge 3^{\alpha}/2$. Combining the upper and the
lower bound for $n$ yields
$3^{\alpha}\le 2\alpha-1$, which has no solution with $\alpha\ge 1$. \qed\\

\noindent {\tt Remark}. It is not difficult to show directly
that if $3\nmid m$, then $2{\rm ord}_9(4m)\le m$ and thus a proof of Lemma \ref{notdrie} can
be given that is free of periodicity considerations and only involves material
from Section \ref{preparations}.

\subsection{$D_S(n)$ is a prime-power}
Assume $9\nmid d$. By Corollary \ref{9} and Corollary \ref{bertie} we can take $g(n)=2n-1$ in
 Lemma \ref{gee}. This yields Lemma \ref{nonnie}. 
However, for the convenience of the reader we give a more direct proof. 
\begin{Lem}\label{nonnie}
Suppose that $d$ with $9\nmid d$ satisfies $\rho(d)\le d/2$, 
then $d$ is a Salajan non-value.
\end{Lem}
\textit{Proof}. Suppose that $d=D_S(n)$ for some integer $n$. By Corollary \ref{bertie} we have $d<2n$.
By Lemma \ref{generalperiod} the condition $9\nmid d$ guarantees that the Salajan sequence is purely
periodic modulo $d$. Since $u_1\equiv u_{1+\rho(d)}({\rm mod~d})$ we must have
$\rho(d)\ge n$. Now suppose that $d\ge 2\rho(d)$. It then follows that $d\ge 2n$, contradicting
$d=D_S(n)<2n$. \qed\\

\indent We now have the necessary ingredients to establish the following result. 
Let $p$ be odd. On noting that in $(\mathbb Z/p^m\mathbb Z)^*$ a square has maximal order $\varphi(p^m)/2$, we
see that the following result says that a Salajan value is either a power of two or 
prime power $p^m$ with $9$ having maximal multiplicative order in $(\mathbb Z/p^m\mathbb Z)^*$.
\begin{Lem}
\label{redpp}
A Salajan value $>1$ must be of the form $p^m$, with $p=2$ or $p>3$ 
and $m\ge 1$. Further, one must have 
ord$_9(p^m)=\varphi(p^m)/2$ and ord$_9(p)=(p-1)/2$. If $m\ge 2$ we
must have $3^{p-1}\not\equiv 1({\rm mod~}p^2)$. 
\end{Lem}
{\it Proof}. Suppose that $d>1$ is a Salajan value that is not a prime power. 
Thus we can write $d=d_1d_2$ with
$d_1,d_2>1$, $(d_1,d_2)=1$. By Lemma \ref{notdrie} we have  $3\nmid d_1d_2$. 
By Lemma \ref{uppierho} we have $\rho(d_1d_2)\le d_1d_2/2$, which by
Lemma \ref{nonnie} implies that $d=d_1d_2$ is a non-value. Thus $d$ is a 
prime power $p^m$. By Lemma \ref{notdrie} we have $p=2$ or $p>3$. Now
let us assume that $p>3$. 
By Lemma \ref{previous} we have either $\rho(p^m)=\varphi(p^m)$ or $\rho(p^m)\le \varphi(p^m)/2$.
The latter inequality leads to $\rho(p^m)\le p^m/2$ and hence to $p^m$ being a non-value.
Using Theorem \ref{generalperiod} we infer that ord$_9(p^m)=\varphi(p^m)/2$. Now if ord$_9(p^m)<(p-1)/2$, this leads
to ord$_9(p^m)<\varphi(p^m)$ and hence we must have ord$_9(p^m)=(p-1)/2$. 
Finally, suppose that $m\ge 2$ and $3^{p-1}\equiv 1({\rm mod~}p^2)$. 
Then ord$_9(p^m)<\varphi(p^m)/2$. This contradiction shows that if $m\ge 2$ we
must have $3^{p-1}\not\equiv 1({\rm mod~}p^2)$. \qed\\

\noindent The possible Salajan values can be further limited by using some results on a quantity
we will baptise as the {\it incongruence index}.

\subsection{$D_S(n)$ is a prime or a small prime power}
\label{primesets}
Put ${\cal P}=\{p:p>3,~{\rm ord}_9(p)=(p-1)/2\}$. 
If a prime $p>3$ is a Salajan value, then by Lemma \ref{redpp} we must have $p\in {\cal P}$.
If $p\in {\cal P}$, then by Theorem \ref{generalperiod} we have $\rho(p)=p-1$. This will be used a few times in the
sequel.
Let 
$${\cal P}_j=\{p:p>3,~p\equiv j({\rm mod~}4),~{\rm ord}_3(p)=p-1\},~j\in \{1,3\}$$ 
and $${\cal P}_2=\{p:p>3,~p\equiv 3({\rm mod~}4),~{\rm ord}_3(p)=(p-1)/2\}.$$
By equation (\ref{lcmorder}) we have $2{\rm ord}_9(p)={\rm lcm}(2,{\rm ord}_3(p))$. {}From
this we infer that ${\cal P}={\cal P}_1 \cup {\cal P}_2 \cup {\cal P}_3$.
We have
$${\cal P}_1=\{5,17,29,53,89,101,113,137,149,173,197,233,257,269,281,293,\ldots\},$$
$${\cal P}_2=\{11,23,47,59,71,83,107,131,167,179,191,227,239,251,263,\ldots\},$$
$${\cal P}_3=\{7,19,31,43,79,127,139,163,199,211,223,283,
\ldots\}.$$
(The reader interested in knowing the natural densities of these sets, under GRH, is referred
to the appendix.)\\
\indent The aim of this section is to establish
the following result, the proof of which makes use of properties of the incongruence index and is
given in Section \ref{5.4.1}.
\begin{Prop}
\label{reducetoprime}
Let $d>1$ be an integer coprime to $10$. If $d$ is a Salajan value, then $d\in {\cal P}_1\cup {\cal P}_2$.
\end{Prop}
\subsubsection{The incongruence index}
\begin{Def}
Let $\{v_j\}_{j=1}^{\infty}$ be a sequence of integers and $m$ an integer. Then 
the largest number $k$ such that $v_1,\ldots,v_{k}$ are pairwise incongruent modulo $m$, we
call the incongruence index, $\iota_v(m)$, of $v$ modulo $m$.
\end{Def}
Note that $\iota_v(m)\le m$. In case the sequence $v$ is purely periodic
modulo $d$, we have $\iota_v(d)\le \rho_v(d)$. A minor change in the proof
of Lemma \ref{gee} yields the following result.
\begin{Lem}
\label{gee2}
Assume that $D_v(n)\le g(n)$ for every $n\ge 1$ with $g$ non-decreasing.
If $d>g(\iota_v(d))$, then $d$ is a $D_v$-non-value.
\end{Lem}
Likewise a minor variation in the proof of Lemma \ref{nonnie} gives the following result, which will be of vital importance in order to discard possible Salajan values.
(For the Salajan sequence $u$ we write $\iota(d)$ instead of $\iota_u(d)$.) 
\begin{Lem}
\label{nonnie2}
If $\iota(d)\le d/2$, then $d$ is a Salajan non-value.
\end{Lem}

\subsubsection{Lifting from $p^m$ to $p^{m+1}$}

\begin{Lem} 
\label{indelift}
If $p>3$ and $\iota(p^m)<\rho(p^m)$, then $\iota(p^{m+1})<p^{m+1}/2$.
\end{Lem}
{\it Proof}. Either $\rho(p^{m+1})=\rho(p^m)$ or
$\rho(p^{m+1})=p\rho(p^m)$. In the first case
$$\iota(p^{m+1})\le \rho(p^{m+1})=\rho(p^m)\le p^m<p^{m+1}/2,$$ so we may
assume that $\rho(p^{m+1})=p\rho(p^m)$. This implies
that
\begin{equation}
\label{nonmirimanoff}
3^{\rho(p^m)}\equiv 1+kp^m({\rm mod~}p^{m+1})
\end{equation}
with $p\nmid k$. {}From this we infer that
$u_{i+j\rho(p^m)}$ assumes $p$ different values modulo $p^{m+1}$ as
$j$ runs through $0,1,\ldots,p-1$. 
Put $j_1=\iota(p^{m})+1$. By assumption there
exists $1\le i_1<j_1$ such that $u_{i_1}\equiv u_{j_1}({\rm mod~}p^{m})$.
Modulo $p^{m+1}$ we have
$$\{u_{i_1+j\rho(p^m)}:0\le j\le p-1\}=\{u_{j_1+j\rho(p^m)}:0\le j\le p-1\}.$$ 
The cardinality of these sets is $p$. Now let us consider the subsets obtained from the
above two sets if we restrict $j$ to be $\le p/2$. Each contains
$(p+1)/2$ different elements. It follows that these sets must have an element
in common. Say we have
$$u_{i_1+k_1\rho(p^m)}\equiv u_{j_1+k_2\rho(p^m)}({\rm mod~}p^{m+1}),~0\le k_1,k_2\le p/2.$$
Since by assumption $i_1\not\equiv j_1({\rm mod~}\rho(p^m))$, we have that
$$i_1+k_1 \rho(p^m)\ne j_1+k_2\rho(p^m).$$
The proof is completed on noting that 
$i_1+k_1 \rho(p^m)$ and $j_1+k_2\rho(p^m)$ are bounded above by
$$\iota(p^m)+1+(p-1){\rho(p^m)\over 2}\le (p+1){\rho(p^m)\over 2}\le (p+1){\varphi(p^m)\over 2}=
p^{m-1}{(p^2-1)\over 2}<{p^{m+1}\over 2},$$
where we used 
that by assumption $\iota(p^m)+1\le \rho(p^m)$ and Lemma \ref{previous}. \qed 

\begin{Lem}
\label{Bb}
Suppose that $l\ge 1$. If $\iota(p^l)\le p^l/2$, then $\iota(p^m)\le p^m/2$ for every $m>l$.
\end{Lem}
{\it Proof}. Note that $p>5$. If $p\not\in {\cal P}$, then $\rho(p^{m})\le p^{m-1}\rho(p)\le p^{m-1}(p-1)/2
\le p^m/2$ and hence  $\iota(p^m)\le \rho(p^m)\le p^m/2$, so we may assume that $p\in {\cal P}$.
Now we proceed by induction. Suppose that we have established that $\iota(p^k)\le p^k/2$ for
$l\le k\le m-1$. By Corollary \ref{liftje2} there are two cases to be considered.\\
Case 1. $\rho(p^2)=\rho(p)=p-1$.\\
In this case $\rho(p^m)\le p^{m-2}\rho(p)=\varphi(p^{m-1})\le p^{m-1}\le p^m/2,$
and hence $\iota(p^{m})\le p^m/2$.\\
Case 2. We have $\rho(p^2)=p\rho(p)$ and hence $\rho(p^{m})=p^{m-1}\rho(p)=\varphi(p^{m})$.
By assumption we have $\iota(p^{m-1})\le p^{m-1}/2<p^{m-2}(1-1/p)=\rho(p^{m-1})$.
By Lemma \ref{indelift} it then follows that $\iota(p^m)\le p^m/2$. \qed\\

\noindent On combining the latter two lemmas with Lemma \ref{nonnie2} we arrive at the following 
more appealing result.
\begin{Lem}\label{18}$~$\\
\noindent {\rm 1)} If $p>3$ and $\iota(p)<\rho(p)$, then $p^2,p^3,\ldots$ are all Salajan non-values.\\
{\rm 2)} If $\iota(p)\le p/2$, then $p,p^2,p^3,\ldots$ are all Salajan non-values.
\end{Lem}
{\it Proof}. 1) If the conditions on $p$ are satisfied, then 
by Lemma \ref{indelift} it follows that $\iota(p^2)\le p^2/2$, 
which by Lemma \ref{Bb} implies that $\iota(p^m)\le p^m/2$ for every $m\ge 2$. 
By Lemma \ref{nonnie2} it then follows that $p^m$ is a non-value.\\
2) If $\iota(p)\le p/2$, then $\iota(p^m)\le p^m/2$ for every $m\ge 1$ by Lemma \ref{Bb} and by 
Lemma \ref{nonnie2} it then follows that $p^m$ is a non-value. \qed\\

\noindent We will see in Proposition \ref{sharpiotap} that actually $\iota(p)\le p/2$ for $p>5$.

\subsection{If ord$_9(p)=(p-1)/2$, then $\iota(p)<\rho(p)$ unless $p=5$}
\label{tip}
Lemma \ref{nonnie2} in combination with the following lemma shows that every $p\in {\cal P}_3$ is a 
Salajan non-value. Recall that if $p\in {\cal P}$, then $\rho(p)=p-1$.
\begin{Lem}
\label{halfpforP_3}
Suppose that $p\in {\cal P}_3$. Then $\iota(p)\le p/2<p-1=\rho(p)$.
\end{Lem}
{\it Proof}. Since by assumption $3$ is a primitive root modulo $p$, 
we have that $({3\over p})=-1$. It follows
that
$$1=u_2={({3\over p})+5\over 4}\equiv {3^{(p-1)/2}+5\over 4}=u_{p-1\over 2}({\rm mod~}p).$$
We infer that $\iota(p)\le (p-1)/2$.\qed

On using Lemma \ref{indelift} the following result can be used to show that
if $p\in {\cal P}_1$ and $m\ge 2$, then $p^m$ is a Salajan non-value.
\begin{Lem}
\label{remainder}
If $p>5$ and $p\in {\cal P}_1\cup {\cal P}_2$, then there exists
$k\le p-3$ such that $u_{k}\equiv u_{k+1}({\rm mod~}p)$ and hence $\iota(p)<p-1=\rho(p)$.
\end{Lem}
{\it Proof}. Note that $$u_{2m-1}\equiv u_{2m} ({\rm mod~}p){\rm ~iff~}3^{2m}\equiv 15({\rm mod~}p)$$ and
$$u_{2m}\equiv u_{2m+1} ({\rm mod~}p){\rm ~iff~}3^{2m}\equiv -5({\rm mod~}p).$$ 
If $p\in {\cal P}_1$, then $3$ is a primitive root modulo $p$, hence $({3\over p})=-1$ and
$({-3\over p})=-1$ as $p\equiv 1({\rm mod~}4)$. If $p\in {\cal P}_2$, then 
$({3\over p})=1$ and $({-3\over p})=-1$ as $p\equiv 3({\rm mod~}4)$.
We see that
$({15\over p})=({-3\over p})({-5\over p})=-({-5\over p})$ and hence 
either $15$ or $-5$ is a square modulo $p$. Since by assumption
ord$_9(p)=(p-1)/2$, every square $s\ne 0$ 
modulo $p$ is of the form $s=3^{2k}$ for some $1\le k\le (p-1)/2$.
It follows that either $3^{2k}\equiv -5({\rm mod~}p)$ or $3^{2k}\equiv 15({\rm mod~}p)$
for some $1\le k\le (p-1)/2$.
Since $3^{p-1}\equiv 1 ({\rm mod~}p)$ and, modulo $p$, $-5$ and $15$ are not congruent to $1$, it
follows that $2k\le p-3$ and so $\iota(p)\le p-3+1=p-2$. 
\qed\\

\noindent {\tt Remark.} We have $({15\over p})=({-5\over p})$ in case $p\in {\cal P}_3$ and
$({-5\over p})=-1$ iff $p\equiv \pm 1({\rm mod~}5)$. We infer that 
if $p>5$ and $p\in {\cal P}$, then there exists
$k\le p-3$ such that $u_{k}\equiv u_{k+1}({\rm mod~}p)$, except when $p\in  {\cal P}_3$ and
$p\equiv \pm 1({\rm mod~}5)$.\\

\noindent {\tt Remark.} It is not true in general that $\iota(p)<\rho(p)$, there are
many counter-examples, e.g., $p=193,307,1093,1181,1871$. It is an open problem whether there
are infinitely many prime numbers $p$ such that $\iota(p)=\rho(p)$.

\subsubsection{Proof of Proposition \ref{reducetoprime}}
\label{5.4.1}
Suppose that $(d,10)=1$. By Lemma \ref{redpp} it follows that $d=p^m$ with $p>5$ and $p\in {\cal P}$.
It follows from Lemmas~\ref{halfpforP_3} and~\ref{remainder} that $\iota(p)<\rho(p)$ for 
every $p\in \mathcal P$ with $p>5$, which implies by Lemma \ref{18} that $m=1$ and $d=p$.

By Lemma \ref{nonnie} and Lemma \ref{halfpforP_3} every prime $p\in {\cal P}_3$ is a
Salajan non-value. On recalling that ${\cal P}={\cal P}_1 \cup {\cal P}_2 \cup {\cal P}_3$ 
the proof is then completed. 
%
\qed

\subsection{$D_S(n)$ is not a `big' prime}  
We will now use classical exponential sum techniques to show that, for sufficiently large primes, the condition given in Corollary~\ref{prop2n} is not satisfied. Therefore, big primes are Salajan non-values.

Let us denote by $\psi$ the additive characters of the group $G$ and $\psi_0$ the trivial character. For any non-empty subset $A\subseteq G$, let us define the quantity \begin{equation}\label{def_S(A)}
|\widehat{A}|=\max_{\psi\neq \psi_0} \left| \sum_{a\in A} \psi(a)\right|,
\end{equation}
where the maximum is taken over all non-trivial characters in $G$.

\begin{Lem}\label{B+B_zero} Let $G$ be a finite abelian group. For any given non-empty subsets $A,B\subseteq G$, whenever $A\cap (B+B)=\emptyset$ we have 
$$|B|\le {|\widehat{A}||G|\over |A|+|\widehat{A}|},$$ where $|\widehat{A}|$ is the quantity defined in~\eqref{def_S(A)}.
\end{Lem}
{\it Proof}. The number $N$ of pairs $(b,b')\in B\times B$ such that $b+b'\in A$ equals
\begin{equation}
\label{characters}
N=\frac{1}{|G|}\sum_{\psi}\sum_A\sum_{B\times B} \psi(b+b' -a) = \frac{|B|^2|A|}{|G|} + R
\end{equation}
where, by the orthogonality of the characters,
\begin{align*}
|R|& = \left| \frac{1}{|G|}\sum_{\psi\ne \psi_0}\sum_A\sum_{B\times B}
\psi(b+b' -a) \right|
\le \frac{1}{|G|}\sum_{\psi\ne \psi_0}\left|\sum_A
\psi(a)\right|\left|\sum_B\psi(b)\right|^2 \\
&\le \frac{|\widehat{A}|}{|G|}\sum_{\psi\ne
\psi_0}\left|\sum_B\psi(b)\right|^2.
\end{align*}
Note that
\[
\left| \sum_B \psi(b) \right|^2 = \sum_{b,b'\in B}\psi (b-b'),
\]
since as complex numbers $\overline{\psi(b)}=\psi(-b)$, and that by
orthogonality of the characters
\[
\sum_{\psi} \sum_{b,b'\in B} \psi(b-b') = \left\{
\begin{array}{lc} 0&\text{ if } b\neq b',\\
|G|& \text{ if } b=b'.
\end{array}
\right.
\]
Thus
\begin{equation}
\label{Error}   
|R|\le \frac{|\widehat{A}|}{|G|} \sum_{\psi\neq \psi_0} \left| \sum_B
\psi(b) \right|^2 = \frac{|\widehat{A}|}{|G|} \left( |G||B| -
|B|^2\right).
\end{equation}
Since by assumption $N=0$, it follows from~\eqref{characters} and~\eqref{Error} that 
\[\frac{|B|^2|A|}{|G|}\le \frac{|\widehat{A}|}{|G|}(|G||B|-|B|^2),\]
which concludes the proof.
\qed

We will need the following auxiliary result, which can be found in~\cite{Ana}.
\begin{Lem}\label{S(A)} Let $p$ be a prime and $g$ be a primitive root modulo $p$. The set
\[
A=\{ (x,y):\, 3g^x-g^{y} \equiv 30 \pmod p\}\subset \mathbb Z_{p-1} \times \mathbb Z_{p-1}
\]
has $p-2$ elements and satisfies $|\widehat{A}|<p^{1/2}$.
\end{Lem} 

\noindent {\tt Remark}. It is easy to see that any subset of an abelian group satisfies that $|A|^{1/2}\le |\widehat{A}|$, so the bound in Lemma~\ref{S(A)} is essentially best possible.

\smallskip

\noindent {\tt Remark}. In fact this result is true in a more general context (see for example~\cite{Ana}): let $g$ be a primitive root in a finite field $\mathbb F_q$ and $a$, $b$ and $c$ be non-zero elements in the field. Then, the set $A_g(a,b,c)=\{(x,y): ag^x-bg^y = c\}$ in $\mathbb F_q$ has $q-2$ elements and satisfies $|\widehat{A}_g(a,b,c)|<q^{1/2}$.

\begin{Prop}\label{prime_lemma} Let $p>3$ be a prime. 
Suppose that $u_1,\ldots,u_n$ are pairwise distinct modulo $p$. Then $p> \left\lfloor \frac{n}{4} \right\rfloor^{4/3}$.
\end{Prop}
\noindent {\it Proof}. First observe that if two elements have the same
parity index, then
$u_i\not\equiv u_{i+2k} ({\rm mod~}p)$ iff $
9^k\not\equiv 1({\rm mod~}p),$ thus ord$_9(p)\ge n/2$. 
(Alternatively one might invoke Lemma \ref{oneventje} to obtain this conclusion.)
By hypothesis,
comparing elements with distinct parity index, it follows that
\begin{equation}
\label{power3} 3\cdot 9^{k}-9^{s} \equiv 30 \pmod p,\ 1\leq
k,s\leq \left\lfloor \tfrac{n}{2}\right\rfloor
\end{equation}
has no solution (otherwise $u_{2k}\equiv u_{2s-1} \pmod p$, with $1\le 2k,
2s-1\le n$).

We will now show that the non existence of solutions to
equation~\eqref{power3} implies that $p >\lfloor \frac{n}{4} \rfloor^{4/3}$.
Let $g$ be a primitive root modulo $p$ and let $A$ be the set defined in
Lemma~\ref{S(A)}. Let $m$ be the smallest integer such that
$g^m\equiv 9\pmod p$ and 
$$B=\{(mx,my):1\le x,y\le  \lfloor n/4\rfloor\} \subset \mathbb Z_{p-1} \times \mathbb Z_{p-1}.$$
Note that, since ord$_9(p)\ge n/2$, it follows that $|B|=\left\lfloor
\frac{n}{4} \right\rfloor^2$ (since $m$ generates a subgroup of order at
least $n/2$ modulo $p-1$).

Observe that the non existence of solutions to
equation~\eqref{power3} implies that
\begin{equation*}
\label{power4} 3\cdot g^{mk}-g^{ms} \equiv 30 \pmod p,\ 1\leq
k,s\leq \left\lfloor \tfrac{n}{2}\right\rfloor
\end{equation*}
has no solutions and in particular $A\cap (B+B)=\emptyset$ (since clearly
$B+B\subseteq \{(mx,my): 1\le x,y\le \left\lfloor n/2\right\rfloor \}$). It
follows from Lemma~\ref{B+B_zero} and Lemma~\ref{S(A)} that
\begin{equation}
\label{bijna}
|B|=\left\lfloor \frac{n}{4} \right\rfloor^2 \le
\frac{|\widehat{A}||G|}{|A|+|\widehat{A}|}\le
\frac{p^{1/2}(p-1)^2}{p-2+p^{1/2}}<
p^{3/2},
\end{equation}
which concludes the proof.
\qed

\begin{cor}\label{Nobigprimes}
If $p > 5$ is a prime number, then $p$ is a Salajan non-value.
\end{cor}
\textit{Proof.} First observe that, if $n\ge 2060$ then it follows from Proposition~\ref{prime_lemma}  that if, for some prime $p\ge n$ the elements $u_1,\ldots,u_n$ are pairwise distinct modulo $p$ then
\[
p> \left\lfloor \frac{n}{4} \right\rfloor^{4/3} \ge 2n,
\]
and by Corollary~\ref{prop2n} it follows that $p$ is a Salajan non-value. For primes $5\le p \le 2060$, the result follows from the calculations included in Table 1.
\qed\\

\noindent Taking $n=\iota(p)$ in Proposition \ref{prime_lemma} we obtain, after
some numerical work, the following estimate. Since $\iota(29)=14$ the 
bound is sharp.
\begin{Prop}
\label{sharpiotap}
Let $p>5$ be a prime. Then
$\iota(p)\le \min((p-1)/2,4p^{3/4})$.
\end{Prop}
{\it Proof}. By Proposition \ref{prime_lemma} we infer that $\iota(p)<3+4p^{3/4}$. A tedious analysis using 
the one but last estimate for $|B|$ in (\ref{bijna}) gives the more elegant bound $\iota(p)<4p^{3/4}$.
For $p<4111$ one verifies the claimed bound by direct computation.
Since $4p^{3/4}<(p-1)/2$ for $p\ge 4111$, we are done. \qed

\section{The proof of Salajan's conjecture}
In Section~\ref{preparations}, we established that powers of $2$ and powers of $5$ were candidates for Salajan values. Finally, after studying the characteristics of the period and the incongruence index of the Salajan sequence, we discard in Section~\ref{non-values} any other possible candidates.

\medskip

\noindent {\it Proof of Theorem} \ref{main}.
It follows from Proposition~\ref{reducetoprime} that if $d>1$ is a Salajan value, then either $(10,d)>1$ or $d\in \mathcal P_1\cup \mathcal P_2$. It follows from Corollary~\ref{Nobigprimes} that no prime greater than $5$ can be a Salajan 
value and hence $(10,d)>1$. By Lemma~\ref{nonnie} it follows that $d$ has to be a prime power.
Therefore, since $(10,d)>1$, the discriminator must be a power of $2$ or a power of $5$.

First suppose that $D_S(n)=2^e$. On invoking Lemma \ref{prop2n} it then follows
that $e=\min\{a:2^a\ge n\}$. Next suppose that $D_S(n)=5^f$. By Lemma
\ref{prop5n} it then follows that $f=\min\{a:2^a\ge 5n/4\}$. So we have
$D_S(n)=2^e$ or $D_S(n)=5^f$. By the definition of the discriminator
we now infer that $D_S(n)=\min\{2^e,5^f\}$. \qed

\section{Appendix}
\subsection{The natural density of the sets ${\cal P}_i$}
Standard methods allow one to determine, assuming the Generalized Riemann Hypothesis, 
the densities of the sets ${\cal P}_i$ defined 
in Section \ref{primesets}. (For a survey of related material see Moree \cite{Asurvey}.)
\begin{Prop}
Assume GRH. We have
$$\#\{p\le x:p\in {\cal P}_i\}=\delta({\cal P}_i){x\over \log x}+O\Big({x\log \log x\over \log^2 x}\Big),$$
with $\delta({\cal P}_1)=\delta({\cal P}_2)=3A/5=0.224373488\ldots$ and 
$\delta({\cal P}_3)=2A/5=0.149582325\ldots$ and
$$A=\prod_p\left(1-{1\over p(p-1)}\right)=0.3739558136\ldots,$$
the Artin constant.
\end{Prop}
\begin{cor}
 The result also holds for the set ${\cal P}$, where 
we find $\delta({\cal P})=\delta({\cal P}_1)+\delta({\cal P}_2)+\delta({\cal P}_3)=8A/5=0.598329301\ldots$.
\end{cor}
{\it Proof}. These three results can be obtained by a variation of the classical result of
Hooley \cite{H} and this yields the estimate with $\delta({\cal P}_i)$ yet to be determined. We
note that the sets ${\cal P}_i$ are mutually disjunct. 
By \cite[Theorem 4]{PU} we have $\delta({\cal P}_1)=3A/5$ and $\delta({\cal P}_1\cup {\cal P}_3)=A$.
This gives $\delta({\cal P}_3)=2A/5$. By \cite[Theorem 3]{PNear} we have
$\delta({\cal P})=8A/5$ and hence $\delta({\cal P}_2)=\delta({\cal P})-
\delta({\cal P}_1\cup {\cal P}_3)=3A/5$. \qed\\

\indent For the benefit of the reader we give a perhaps more insightful argument why $\delta({\cal P}_2)=3A/5$.\\
\indent Assuming GRH we have, cf. Moree \cite{PNear}, 
$$\delta({\cal P}_2)=\sum_{n=1}^{\infty}{\mu(n)\over [\mathbb Q(\zeta_{2n},3^{1/2n}):\mathbb Q]}
-\sum_{n=1}^{\infty}{\mu(n)\over [\mathbb Q(i,\zeta_{2n},3^{1/2n}):\mathbb Q]},$$
where the first sum gives the density of the primes $p$ such that ord$_3(p)=(p-1)/2$ and the second
sum the density of the primes $p$ such that $p\equiv 1({\rm mod~}4)$ and ord$_3(p)=(p-1)/2$.
Since for $n$ even, $i\in \mathbb Q(\zeta_{2n})$, we find that
$$\delta({\cal P}_2)=\sum_{(n,2)=1}^{\infty}{\mu(n)\over [\mathbb Q(\zeta_{2n},3^{1/2n}):\mathbb Q]}
-\sum_{(n,2)=1}^{\infty}{\mu(n)\over [\mathbb Q(i,\zeta_{2n},3^{1/2n}):\mathbb Q]}.$$
Now suppose that $n$ is odd. If $3|n$, then $\sqrt{-3}\in \mathbb Q(\zeta_{2n})$. Since
$\sqrt{3}\in \mathbb Q(\zeta_{2n},3^{1/2n})$, it follows that 
$\mathbb Q(i,\zeta_{2n},3^{1/2n})=\mathbb Q(\zeta_{2n},3^{1/2n})$. On the other hand,
if $(n,3)=1$ one infers that
$[\mathbb Q(i,\zeta_{2n},3^{1/2n}):\mathbb Q]=2[\mathbb Q(\zeta_{2n},3^{1/2n}):\mathbb Q]$.
This leads to
$$\delta({\cal P}_2)={1\over 2}\sum_{(n,6)=1}^{\infty}{\mu(n)\over [\mathbb Q(\zeta_{2n},3^{1/2n}):\mathbb Q]}
={1\over 4}\sum_{(n,6)=1}{\mu(n)\over n\varphi(n)}={3\over 5}A,$$
where we used that $[\mathbb Q(\zeta_{2n},3^{1/2n}):\mathbb Q]=\varphi(2n)2n=2\varphi(n)n$ if $(n,6)=1$
and the identity
$$\sum_{(n,6)=1}{\mu(n)\over n\varphi(n)}=\prod_{p>3}\left(1-{1\over p(p-1)}\right)={12\over 5}A.$$

\subsection{Counting the elements $\le x$ in ${\cal F}$}
\label{Izabela}
In this section, written jointly with Izabela Petrykiewicz, we will establish 
Proposition \ref{iza} from the introduction.\\
\indent Recall that 
${\cal F}=\{f~:~[4\cdot 5^{f-1}, 5^f]\textnormal{ contains no power of 2}\}$. 
Consider ${\cal G}={\mathbb N}\setminus {\cal F}$. 
We have that $g$ is in ${\cal G}$ iff $4\cdot 5^{g-1}\leq 2^k \leq 5^g$ for some $k\in {\mathbb N}$.
Thus we have $g$ is in ${\cal G}$ iff $2\log2 + (g-1)\log 5 \leq k \log 2 \leq g\log 5$, that
is iff
$2+(g-1) \alpha\leq k \leq g \alpha,$
where $\alpha=\log 5/\log 2$. Since $k$ is an integer, we may replace $g\alpha$ by $[g\alpha]$ and
the condition becomes $k\in [[g\alpha]+\{g\alpha\}+2-\alpha,[g\alpha]]$. Note that there can be only
an integer in this interval iff $\{g\alpha\}\le \alpha-2$.
Note that $\alpha$ is irrational. Now it is a consequence of Weyl's criterion, see, e.g., \cite{Bug,KN}, that 
for a fixed $0<\beta<1$ we have
$$\#\{g\le x:\{g\alpha\}\le \beta\}\sim \beta x,~x\rightarrow \infty.$$ On applying this with $\beta=\alpha-2$ the proof of Proposition \ref{iza} is 
easily completed. \qed\\

\noindent {\tt Acknowledgement}. This project was started in the context
of an internship of Sabin Salajan at MPIM in 2012 and a visit of the
first author to the second author to ICMAT in Madrid. The project was
taken up again during a two week visit of Bernadette Faye (Senegal) 
at MPIM. The first author thanks Bernadette for discussions and her
help with some computer experiments. Further he thanks David Brink, Igor
Shparlinski and
Arne Winterhof for helpful e-mail correspondence and Paul Tegelaar for comments on
an earlier version. The idea of the proof
of Proposition \ref{iza} is due to Izabela Petrykiewicz.


\begin{thebibliography}{99}
\bibitem{Apostel} T. Apostol, {\it Introduction to analytic number theory}, 
Undergraduate Texts in Mathematics, Springer-Verlag, New York-Heidelberg, 1976.
\bibitem{ABM} L.K. Arnold, S.J. Benkoski and B.J. McCabe, The
discriminator (a simple application of Bertrand's postulate),
{\it Amer. Math. Monthly} {\bf 92} (1985), 275–-277.
\bibitem{Beyl} R.F. Beyl, Cyclic subgroups of the prime residue group,
{\it Amer. Math. Monthly}
{\bf 84} (1977), 46–-48.
\bibitem{BSW} P.S. Bremser, P.D. Schumer and L.C. Washington, A note
on the incongruence of consecutive
integers to a fixed power, {\it J. Number Theory} {\bf 35} (1990),
105–-108.
\bibitem{BC} J. Browkin and H.-Q. Cao, Modifications of the Eratosthenes sieve, 
{\it Colloq. Math.} {\bf 135} (2014), 127--138. 
\bibitem{Bug} Y. Bugeaud, {\it Distribution modulo one and Diophantine approximation}, Cambridge Tracts in Mathematics {\bf 193}, Cambridge University Press, Cambridge, 2012. 
\bibitem{Ana} J. Cilleruelo and A. Zumalac\'arregui, An additive problem in finite fields with powers 
of elements of large multiplicative order, {\it Rev. Mat. Complut.} {\bf 27} (2014), 501--508.
\bibitem{H} C. Hooley, On Artin's conjecture, 
{\it J. Reine Angew. Math.} {\bf 225} (1967), 209--220. 
\bibitem{KN} L. Kuipers and H. Niederreiter, {\it Uniform distribution of sequences}, Pure and Applied Mathematics, 
Wiley-Interscience, New York-London-Sydney, 1974
\bibitem{Mann} H.B. Mann, {\it Addition theorems: The addition
theorems of group theory and number theory},
Interscience Publishers John Wiley and Sons, New York-London-Sydney,
1965.
\bibitem{PP} P. Moree, The incongruence of consecutive values of
polynomials, {\it Finite Fields Appl.}
{\bf 2} (1996), 321–-335.
\bibitem{PU} P. Moree, Uniform distribution of primes having a prescribed primitive root, 
{\it Acta Arith.} {\bf 89} (1999), 9--21. 
\bibitem{Asurvey} P. Moree, Artin's primitive root conjecture -- a survey, 
{\it Integers} {\bf 12A} (2012), No. 6, 1305--1416.
\bibitem{PNear} P. Moree, Near-primitive roots, {\it Funct. Approx. Comment. Math.} {\bf 48.1} (2013), 133--145.
\bibitem{MM} P. Moree and G. L. Mullen, Dickson polynomial
discriminators,
{\it J. Number Theory} {\bf 59} (1996), 88--105.
\bibitem{Sun} Zhi-Wei Sun, On functions taking only prime values,
 {\it J. Number Theory} {\bf 133} (2013), 2794--2812. 
\bibitem{Zieve} M. Zieve, A note on the discriminator, {\it J. Number
Theory} {\bf 73} (1998), 122–-138.
\end{thebibliography}
\end{document}